\documentclass[12pt]{amsart}
\usepackage{amssymb}
\usepackage{latexsym}
\usepackage{amsmath}

\newtheorem{Theorem}{Theorem}[section]

\newtheorem{Corollary}[Theorem]{Corollary}
\newtheorem{Proposition}[Theorem]{Proposition}

\newtheorem{Example}[Theorem]{Example}

\newenvironment{pf}{\medskip\par\noindent{\bf Proof\/}.}{\hfill
$\Box$}

\begin{document}
\title{On Certain Monomial Sequences}
\author{Zhongming Tang}

\begin{abstract}
We give equivalent conditions for a monomial sequence to be a
d-sequence or a proper sequence, and a sufficient condition for a
monomial sequence to be an s-sequence in order to compute
invariants of the symmetric algebra of the ideal generated by it.
\end{abstract}

\thanks{
\indent Supported by the National Natural Science Foundation of
China.}
\date{}
\maketitle

\section{Introduction}
d-Sequences, proper sequences and s-sequences are three kind
sequences related to symmetric algebras. Firstly, let us recall
their definitions. Let $R$ be a commutative Noetherian ring and
$a_1,\ldots,a_n\in R$. We say that $a_1,\ldots,a_n$ is a
d-sequence if $a_1,\ldots,a_n$ is a minimal generating set of the
ideal $(a_1,\ldots,a_n)$ and
$$
(a_1,\ldots,a_i):a_{i+1}a_k=(a_1,\ldots,a_i):a_k,\,
i=0,\ldots,n-1,\, k\geq i+1,
$$
cf., \cite{Hu1}. Proper sequences were introduced in \cite{HSV}.
$a_1,\ldots,a_n$ is called a proper sequence if
$$
a_{i+1}H_j(a_1,\ldots,a_i;R)=0,\, i=0,\ldots,n-1,\, j>0,
$$
where $H_j(a_1,\ldots,a_i;R)$ denotes the $j$-th Koszul homology
of $a_1,\ldots,a_i$. It is shown in \cite{K} that $a_1,\ldots,a_n$
is a proper sequence if and only if
$a_{i+1}H_1(a_1,\ldots,a_i;R)=0$, $i=0,\ldots,n-1$. We will use
this condition as an equivalent definition of proper sequences.
Recently, a new sequence, so-called s-sequence, was introduced in
\cite{HRT} to study symmetric algebras. Let $M=(f_1,\ldots,f_n)$
be a finitely generated $R$-module with relation matrix
$(a_{ij})_{m\times n}$. Then
$$
\mbox{Sym}(M)=R[y_1,\ldots,y_n]/J,
$$
where $J=(g_1,\ldots,g_m)$ and $g_i=\sum_{j=1}^na_{ij}y_j$,
$i=1,\ldots,m$. Let $<$ be a monomial order on the monomials in
$y_1,\ldots,y_n$ with the property $y_1<\cdots<y_n$. Set
$I_i=(f_1,\ldots,f_{i-1}):f_i$. Then
$(I_1y_1,\ldots,I_ny_n)\subseteq \mbox{in}_<(J)$. We call
$f_1,\ldots,f_n$ an s-sequence (with respect to $<$) if
$$
\mbox{in}_<(J)=(I_1y_1,\ldots,I_ny_n).
$$
If, in addition $I_1\subseteq \cdots\subseteq I_n$, then
$f_1,\ldots,f_n$ is called a strong s-sequence.

In this paper, we consider monomial d-sequences, proper sequences
and s-sequences. Throughout let $R=K[x_1,\ldots,x_m]$ be a
polynomial ring, where $K$ is a field.

In the monomial case, s-sequences can be characterized by
Gr\"{o}bner bases. Let $f_1,\ldots,f_n\in R$ be monomials. For any
$i\neq j$, denote the greatest common divisor of $f_i$ and $f_j$
by $[f_i,f_j]$ and set $f_{ij}=\frac{f_i}{[f_i,f_j]}$. Then $J$ is
generated by $g_{ij}:=f_{ij}y_j-f_{ji}y_i$, $1\leq i<j\leq n$, and
$f_1,\ldots,f_n$ is an s-sequence if and only if $g_{ij},\, 1\leq
i<j\leq n$, form a Gr\"{o}bner basis with respect to some monomial
order $<$ such that $y_1<\cdots<y_n$ and $x_i<y_j$ for any $i,j$,
cf.,\cite{HRT}. In the following, we will assume such an order.

In sections 2 and 3, we will characterize monomial d-sequences and
proper sequences. Let $f_1,\ldots,f_n$ be a minimal monomial
sequence. We show that $f_1,\ldots,f_n$ is a d-sequence if and
only if $[f_i,f_j]\mid f_k$, for all $i<j<k$, and
$[f_i,f_j]=[f_i,f_j^2]$, for all $i<j$, while for $f_1,\ldots,f_n$
to be a proper sequence it is equivalent to that $[f_i,f_j]\mid
f_k$, for all $i<j<k$. Section 4 is devoted to discuss monomial
s-sequences. In \cite{HRT}, a first sufficient condition is given.
In this paper we introduce other sufficient conditions. In
particular we study the condition $[f_i,f_j]\mid f_k$, for any
$i<j<k$. For this sequence, very strong properties of the
annihilator ideals follow, in particular to be a strong
s-sequence. We give a sufficient condition for a monomial sequence
to be an s-sequence which is satisfied by strong s-sequences and
show that this condition is also necessary in the case the
sequence is squarefree and has 4 elements. Finally, we deduce some
properties of the symmetric algebra of an ideal generated by a
strong monomial s-sequence.

\section{Monomial d-Sequences}
In this section, we will give a characterization of monomial
d-sequence.

Let $f_1,\ldots,f_n\in R$ be a monomial sequence. We say that
$f_1,\ldots,f_n$ is minimal if it is a minimal generating set  of
the ideal $(f_1,\ldots,f_n)$, which is equivalent to that there is
no $i\not=j$ such that $f_i\mid f_j$.

By definition, $f_1,\ldots,f_n$ is a d-sequence if and only if
$$
(f_1,\ldots,f_i):f_{i+1}f_k=(f_1,\ldots,f_i):f_k,\, i=1,\ldots,
n-1,\,k\geq i+1.
$$
Note that, for any monomial sequence $g_1,\ldots,g_m,g_{m+1}\in
R$,
$(g_1,\ldots,g_m):g_{m+1}=(\frac{g_1}{[g_1,g_{m+1}]},\ldots,\frac{g_m}{[g_m,g_{m+1}]})$.
Then $f_1,\ldots,f_n$ is a d-sequence if and only if
$$
\left(\frac{f_1}{[f_1,f_{i+1}f_k]},\ldots,\frac{f_i}{[f_i,f_{i+1}f_k]}\right)=
\left(\frac{f_1}{[f_1,f_k]},\ldots,\frac{f_i}{[f_i,f_k]}\right)
,\, i=1,\ldots, n-1,\,k\geq i+1.
$$

The sufficient part of the following theorem is observed in
\cite{HSV}, in which it is shown that, when
$[f_i,f_j]=[f_i,f_j^2]$ for all $i<j$, $f_1,\ldots,f_n$ is a
d-sequence if and only if $[f_i,f_j]\mid f_k$ for all $i<j<k$.

\begin{Theorem}
\label{d-seq} Let $f_1,\ldots,f_n$ be a monomial sequence. Then
$f_1,\ldots,f_n$ is a d-sequence if and only if there is no
$i\not=j$ such that $f_i\mid f_j$ and
\begin{eqnarray*}
&& [f_i,f_j]\mid f_k,\, 1\leq i<j<k\leq n,\\
&& [f_i,f_j]=[f_i,f_j^2],\, 1\leq i<j\leq n.
\end{eqnarray*}
\end{Theorem}

\begin{pf}
Note that, for any $j\leq i$, $[f_j,f_{i+1}f_k]=[f_j,f_k]$ holds
if $k=i+1$ and $[f_j,f_k]=[f_j,f_k^2]$, or if $k>i+1$,
$[f_j,f_k^2]=[f_j,f_k]$ and $[f_j,f_{i+1}f_k]\mid [f_j,f_k^2]$
which is satisfied if $[f_j,f_{i+1}]\mid f_k$. Then the sufficient
part follows.

Now we show the necessary part. We use induction on $n$. When
$n=2$, from $(\frac{f_1}{[f_1,f_2^2]})=(\frac{f_1}{[f_1,f_2]})$ we
see immediately that $[f_1,f_2]=[f_1,f_2^2]$. Assume that $n>2$
and the result is true for $n-1$, thus $[f_i,f_j]\mid f_k$ for all
$i<j<k\leq n-1$ and $[f_i,f_j]=[f_i,f_j^2]$ for all $i<j\leq n-1$.
We are going to show that $[f_i,f_j]\mid f_n$ for all $i<j\leq
n-1$ and $[f_i,f_n]=[f_i,f_n^2]$ for all $i\leq n-1$.

We first use induction on $j$ to show that $[f_i,f_j]\mid f_n$ for
all $i<j\leq n-1$. When $j=2$, from $(\frac{f_1}{[f_1,f_2f_n]})=
(\frac{f_1}{[f_1,f_n]})$ we have $[f_1,f_n]=[f_1,f_2f_n]$, but
$[f_1,f_2]\mid [f_1,f_2f_n]$, hence $[f_1,f_2]\mid [f_1,f_n]$, so
$[f_1,f_2]\mid f_n$. Now assume that $j>2$. For any $i<j\leq n-1$,
from
$$
\left(\frac{f_1}{[f_1,f_jf_n]},\ldots,\frac{f_{j-1}}{[f_{j-1},f_jf_n]}\right)=
\left(\frac{f_1}{[f_1,f_n]},\ldots,\frac{f_{j-1}}{[f_{j-1},f_n]}\right),
$$
we have
$$
\frac{f_i}{[f_i,f_jf_n]}\in
\left(\frac{f_1}{[f_1,f_n]},\ldots,\frac{f_{j-1}}{[f_{j-1},f_n]}\right).
$$
Then $\frac{f_k}{[f_k,f_n]}\mid \frac{f_i}{[f_i,f_jf_n]}$ for some
$k\leq j-1$. We claim that $k=i$. If $k\neq i$, then, by induction
hypothesis on $n$, $[f_k,f_i]\mid f_j$, thus $[f_k,f_i]\mid
[f_i,f_jf_n]$. By induction hypothesis on $j$, $[f_k,f_i]\mid
[f_k,f_n]$, hence $\frac{f_k}{[f_k,f_n]}$ and $
\frac{f_i}{[f_i,f_jf_n]}$ are coprime. Then
$\frac{f_k}{[f_k,f_n]}=1$, hence $f_k\mid f_n$, a contradiction.
Hence $\frac{f_i}{[f_i,f_n]}\mid \frac{f_i}{[f_i,f_jf_n]}$, thus
$[f_i,f_jf_n]\mid [f_i,f_n]$, but $[f_i,f_j]\mid [f_i,f_jf_n]$, we
have $[f_i,f_j]\mid f_n$, as required.

Now we show that $[f_i,f_n]=[f_i,f_n^2]$ for all $i\leq n-1$. From
$$
\left(\frac{f_1}{[f_1,f_n^2]},\ldots,\frac{f_{n-1}}{[f_{n-1},f_n^2]}\right)=
\left(\frac{f_1}{[f_1,f_n]},\ldots,\frac{f_{n-1}}{[f_{n-1},f_n]}\right),
$$
we have
$$
\frac{f_i}{[f_i,f_n^2]}\in
\left(\frac{f_1}{[f_1,f_n]},\ldots,\frac{f_{n-1}}{[f_{n-1},f_n]}\right).
$$
Then there exists some $k\leq n-1$ such that
$\frac{f_k}{[f_k,f_n]}\mid \frac{f_i}{[f_i,f_n^2]}$. We claim
again that $k=i$. If $k\neq i$, then, as $[f_k,f_i]\mid f_n$, we
have $[f_k,f_i]\mid [f_k,f_n]$ and also $[f_k,f_i]\mid
[f_i,f_n^2]$, hence $\frac{f_k}{[f_k,f_n]}$ and
$\frac{f_i}{[f_i,f_n^2]}$ are coprime. Then
$\frac{f_k}{[f_k,f_n]}=1$, so $f_k\mid f_n$, a contradiction.
Hence $\frac{f_i}{[f_i,f_n]}\mid \frac{f_i}{[f_i,f_n^2]}$, thus
$[f_i,f_n^2]\mid [f_i,f_n]$, so $[f_i,f_n]=[f_i,f_n^2]$. The proof
is complete.
\end{pf}

We have immediately from \ref{d-seq} the following

\begin{Corollary}
Let $f_1,\ldots,f_n$ be a monomial sequence. If $f_1,\ldots,f_n$
is a d-sequence, then any subsequence of $f_1,\ldots,f_n$ is also
a d-sequence.
\end{Corollary}

\section{Monomial Proper Sequences}
Let $f_1,\ldots,f_n$ be a monomial sequence. Define homomorphisms
$d_1$ and $d_2$ as follows
\begin{eqnarray*}
d_1\,:\,\bigoplus_{i=1}^nRe_i&\longrightarrow&R\\
e_i&\mapsto&f_i, \\
d_2\,:\,\bigoplus_{1\leq i<j\leq
n}Re_{ij}&\longrightarrow&\bigoplus_{i=1}^nRe_i\\
e_{ij}&\mapsto&f_je_i-f_ie_j.
\end{eqnarray*}
Then
$$
H_1(f_1,\ldots,f_n;R)=\mbox{Ker}(d_1)/\mbox{Im}(d_2),
$$
and
$$
\mbox{Ker}(d_1)=(f_{ij}e_j-f_{ji}e_i\,:\,1\leq i<j\leq n).
$$

The following theorem gives a characterization of monomial proper
sequences.

\begin{Theorem}
\label{proper}
Let $f_1,\ldots,f_n$ be a minimal monomial
sequence. Then $f_1,\ldots,f_n$ is a proper sequence if and only
if
$$
[f_i,f_j]\mid f_k,\,\mbox{ for all } 1\leq i<j<k\leq n.
$$
\end{Theorem}

\begin{pf} Let $d_1$ and $d_2$ be defined as above.
Note that, for $r=2,\ldots,n-1$, to show that
$f_{r+1}H_1(f_1,\ldots,f_r;R)=0$ is equivalent to prove that
$$
f_{r+1}(f_{ij}e_j-f_{ji}e_i)\in\mbox{Im}(d_2),\,\mbox{ for all }
 i<j\leq r.
$$

Suppose that $[f_i,f_j]\mid f_k$, for all $i<j<k$. Then, for
$r=2,\ldots,n-1$ and any $i<j\leq r$,
\begin{eqnarray*}
f_{r+1}(f_{ij}e_j-f_{ji}e_i)&=&
\frac{f_{r+1}}{[f_i,f_j]}(f_ie_j-f_je_i)\\
&=& d_2(-\frac{f_{r+1}}{[f_i,f_j]}e_{ij})\in\mbox{Im}(d_2).
\end{eqnarray*}
Thus $f_1,\ldots,f_n$ is a proper sequence.

Conversely, assume that $f_1,\ldots,f_n$ is a proper sequence. We
use induction on $2\leq r\leq n-1$ to show that
$$
[f_i,f_j]\mid f_{r+1},\,\mbox{ for all } i<j\leq r.
$$

Let $r=2$. By $f_3(f_{12}e_2-f_{21}e_1)\in\mbox{Im}(d_2)$, we have
that $f_3(f_{12}e_2-f_{21}e_1)=r(f_2e_1-f_1e_2)$ for some $r\in
R$. Then $f_3f_{12}=-rf_1$, hence $f_3=-r[f_1,f_2]$, so,
$[f_1,f_2]\mid f_3$.

Now assume that $r\geq 3$ and $[f_i,f_j]\mid f_k$ for all
$i<j<k\leq r$, let us show that $[f_i,f_j]\mid f_{r+1}$ for all
$i<j\leq r$. By assumption, for any $i<j\leq r$,
$f_{r+1}(f_{ji}e_i-f_{ij}e_j)\in\mbox{Im}(d_2)$, we see that
$$
f_{r+1}f_{ij}\in(f_1,\ldots,\widehat{f_j},\ldots,f_r).
$$
Note that, if $[f_i,f_j]\nmid f_{r+1}$, then $f_i\nmid
f_{r+1}f_{ij}$, hence
$f_{r+1}f_{ij}\in(f_1,\ldots,\widehat{f_i},\ldots,\widehat{f_j},\ldots,f_r)$.
Since, for any $1\leq i,j\leq r$, $[f_i,f_j]\mid [f_i,f_r]$, we
only need to show that $[f_i,f_r]\mid f_{r+1}$, $i=1,\ldots,r-1$.

Assume the contrary. Then there exists $s\geq 1$ and $1\leq
k_1<\cdots<k_s\leq r-1$ such that
$$
[f_{k_i},f_r]\nmid f_{r+1},\,i=1,\ldots,s,
$$
and
$$
[f_j,f_r]\mid f_{r+1},\,j\not=k_1,\ldots,k_s.
$$

As $[f_{k_i},f_r]\nmid f_{r+1}$, we have that
$f_{r+1}f_{k_i,r}\in(f_1,\ldots,\widehat{f_{k_i}},\ldots,f_{r-1})$.
For any $j\neq k_1,\ldots,k_s$, we claim that $f_j\nmid
f_{r+1}f_{k_i,r}$. If $f_j\mid f_{r+1}f_{k_i,r}$, then, as
$[f_j,f_{k_i}]\mid f_{r+1}$, $f_{j,k_i}\mid
\frac{f_{r+1}}{[f_j,f_{k_i}]}f_{k_i,r}$, but
$[f_{j,k_i},f_{k_i,r}]=1$, hence $f_{j,k_i}\mid
\frac{f_{r+1}}{[f_j,f_{k_i}]}$, thus $f_j\mid f_{r+1}$, a
contradiction. Hence $f_j\nmid f_{r+1}f_{k_i,r}$ for all $j\neq
k_1,\ldots,k_s$. Thus
$$
f_{r+1}f_{k_i,r}\in(f_{k_1},\ldots,\widehat{f_{k_i}},\ldots,f_{k_s}),\,i=1,\ldots,s.
$$

Note that $s\geq 2$. Let $\ell_1=1$. Then there exists
$\ell_2\neq\ell_1$ such that $f_{k_{\ell_2}}\mid
f_{r+1}f_{k_{\ell_1},r}$ and there exists $\ell_3\neq \ell_2$ such
that $f_{k_{\ell_3}}\mid f_{r+1}f_{k_{\ell_2},r}$, here it is
possible that $\ell_3=\ell_1$. Note that $[f_{ k_{\ell_2}},
f_{k_{\ell_3}}]\nmid f_{r+1}$, otherwise, by $f_{k_{\ell_3}}\mid
f_{r+1}f_{k_{\ell_2},r}$ we have that
$f_{k_{\ell_3},k_{\ell_2}}\mid \frac{f_{r+1}}{[f_{ k_{\ell_3}},
f_{k_{\ell_2}}]}f_{k_{\ell_2},r}$, but
$[f_{k_{\ell_3},k_{\ell_2}},f_{k_{\ell_2},r}]=1$, hence
$f_{\ell_3}\mid f_{r+1}$, a contradiction. On the other hand,
since $[f_{k_{\ell_1}}, f_{k_{\ell_2}}]\mid[f_{k_{\ell_1}},f_r]$
and also $[f_{k_{\ell_1}},
f_{k_{\ell_2}}]\mid[f_{k_{\ell_2}},f_r]$, we have that
$f_{k_{\ell_1},r}\mid f_{k_{\ell_1},k_{\ell_2}}$ and
$f_{k_{\ell_2},r}\mid f_{k_{\ell_2},k_{\ell_1}}$. Then
$f_{k_{\ell_2}}\mid f_{r+1}f_{k_{\ell_1},k_{\ell_2}}$ and
$f_{k_{\ell_3}}\mid f_{r+1}f_{k_{\ell_2},k_{\ell_1}}$. Thus
$[f_{k_{\ell_2}},f_{k_{\ell_3}}]\mid
f_{r+1}f_{k_{\ell_1},k_{\ell_2}}$ and also
$[f_{k_{\ell_2}},f_{k_{\ell_3}}]\mid
f_{r+1}f_{k_{\ell_2},k_{\ell_1}}$, but
$[f_{k_{\ell_1},k_{\ell_2}},f_{k_{\ell_2},k_{\ell_1}}]=1$, we get
that $[f_{k_{\ell_2}},f_{k_{\ell_3}}]\mid f_{r+1}$, a
contradiction.
\end{pf}

Let $f_1,\ldots,f_n$ be a monomial d-sequence. Since a d-sequence
is a proper sequence, it follows from \ref{proper} that
$[f_i,f_j]\mid f_k$ for all $i<j<k$. This gives another proof of
$[f_i,f_j]\mid f_k$ in \ref{d-seq}.

From \ref{proper}, we have immediately

\begin{Corollary}
Let $f_1,\ldots,f_n$ be a minimal monomial sequence. If
$f_1,\ldots,f_n$ is a proper sequence, then any subsequence of
$f_1,\ldots,f_n$ is also a proper sequence.
\end{Corollary}

We say a d-sequence (proper sequence) is unconditional if any
permutation of it is also a d-sequence (proper sequence).

\begin{Corollary} Let $f_1,\ldots,f_n$ be a minimal monomial
sequence.
\begin{itemize}
\item[(1)] $f_1,\ldots,f_n$ is an
unconditional proper sequence if and only if there exist monomials
$d,g_1,\ldots,g_n$ such that $f_i=dg_i$, $i=1,\ldots,n$, and
$g_1,\ldots,g_n$ is a regular sequence.
\item[(2)] $f_1,\ldots,f_n$ is an unconditional d-sequence if and only if
there exist monomials $d,g_1,\ldots,g_n$ such that $f_i=dg_i$,
$[d,g_i]=1$, $i=1,\ldots,n$, and $g_1,\ldots,g_n$ is a regular
sequence.
\end{itemize}
\end{Corollary}

\begin{pf} (1) By \ref{proper}, $f_1,\ldots,f_n$ is an
unconditional proper sequence if and only if
$$
[f_i,f_j]=[f_i,f_k]=[f_j,f_k],\, i<j<k.
$$
For any $i<j$, set $d=[f_i,f_j]$ and $g_i=\frac{f_i}{d}$. Then
$[g_i,g_j]=1$ for all $i\neq j$, which is equivalent to that
$g_1,\ldots,g_n$ is a regular sequence.

(2) In virtue of \ref{d-seq} and \ref{proper}, $f_1,\ldots,f_n$ is
an unconditional d-sequence if and only if it is an unconditional
proper sequence and $[f_i,f_j]=[f_i,f_j^2]$ for all $i<j$. Then,
by (1), it is equivalent to that there exist monomials
$d,g_1,\ldots,g_n$ such that $f_i=dg_i$, $i=1,\ldots,n$,
$g_1,\ldots,g_n$ is a regular sequence and
$[dg_i,d^2g_j^2]=[dg_i,dg_j]$, but the latter identity is
equivalent to $[d,g_i]=1$ when $g_1,\ldots,g_n$ is a regular
sequence.
\end{pf}

Let $f_1,\ldots,f_n$ be a squarefree minimal monomial sequence.
Then $[f_i,f_j]=[f_i,f_j^2]$ holds for any $i<j$, thus (1) of
following corollary is clear from \ref{d-seq} and \ref{proper},
but note that a d-sequence is a proper sequence in general. Let
$I=(f_1,\ldots,f_n)$ and $f_i^*$ the canonical image of $f_i$ in
$\mbox{Sym}(I)$. Since $f_1^*,\ldots,f_n^*$ is a d-sequence in
$\mbox{Sym}(I)$ if $f_1,\ldots,f_n$ is a d-sequence by
\cite[Theorem 1.1]{Hu2} and that $f_1^*,\ldots,f_n^*$ is a
d-sequence in $\mbox{Sym}(I)$ is equivalent to that
$f_1,\ldots,f_n$ is a proper sequence by \cite[Theorem 2.2]{K},
the following (2) follows from (1).

\begin{Corollary}
Let $f_1,\ldots,f_n$ be a squarefree minimal monomial sequence.
Then
\begin{itemize}
\item[(1)] $f_1,\ldots,f_n$ is a proper sequence if and only if
$f_1,\ldots,f_n$ is a d-sequence.
\item[(2)] $f_1,\ldots,f_n$ is a d-sequence if and only if
$f_1^*,\ldots,f_n^*$ is a d-sequence in $\mbox{Sym}(I)$, where
$I=(f_1,\ldots,f_n)$.
\end{itemize}
\end{Corollary}

\section{Monomial s-Sequences}
Let $f_1,\ldots,f_n$ be a monomial sequence. Proposition 1.7 of
\cite{HRT} states that $f_1,\ldots,f_n$ is an s-sequence if
$[f_{ij},f_{kl}]=1$ for all $i<j$, $k<l$, $i\not=k$ and $j\neq l$.
The following theorem also gives a sufficient condition for a
monomial sequence to be an s-sequence.

\begin{Theorem}
\label{s-seq} Let $f_1,\ldots,f_n$ be a monomial sequence. For any
$i<j$, $k<l$, $j<l$ and $k\neq i$, suppose that
$[f_{ij},f_{kl}]=1$ or $f_{jl}[f_{ij},f_{kl}]\mid f_{kl}f_{ji}$,
or $f_{ki}[f_{ij},f_{kl}]\mid f_{kl}f_{ji}$ in case $i>k$. Then
$f_1,\ldots,f_n$ is an s-sequence.
\end{Theorem}

\begin{pf}
For any $i<j$ and $k<l$, we need to show that all the S-pairs
$S(g_{ij},g_{kl})$ have standard expressions with zero remainder.

Note that
$$
S(g_{ij},g_{kl})=\frac{f_{kl}f_{ji}}{[f_{ij},f_{kl}]}y_iy_l-
\frac{f_{lk}f_{ij}}{[f_{ij},f_{kl}]}y_ky_j.
$$

If $i=k$, then
$$
S(g_{ij},g_{kl})=y_k(\frac{f_{kl}f_{jk}}{[f_{kj},f_{kl}]}y_l-
\frac{f_{lk}f_{kj}}{[f_{kj},f_{kl}]}y_j).
$$
As $S(g_{ij},g_{kl})=0$ when $y_i=f_i$, $y_j=f_j$, $y_k=f_k$, and
$y_l=f_l$, we see that $\frac{f_{kl}f_{jk}}{[f_{kj},f_{kl}]}y_l-
\frac{f_{lk}f_{kj}}{[f_{kj},f_{kl}]}y_j$ is in the first syzygy of
$f_j,f_l$, thus it is a multiple of $g_{jl}$. Then
$S(g_{ij},g_{kl})$ has a standard expression with zero remainder.
Similarly, $S(g_{ij},g_{kl})$ has a standard expression with zero
remainder if $j=l$.

Now assume that $i<j$, $k<l$, $i\neq k$ and $j\neq l$, and assume
that $j<l$. Then, by assumption,
$\frac{f_{kl}f_{ji}}{[f_{ij},f_{kl}]}$ is divided by $f_{jl}$ or
$f_{kl}$ or $f_{ki}$ when $i>k$. If
$\frac{f_{kl}f_{ji}}{[f_{ij},f_{kl}]}$ is divided by $f_{jl}$,
then
$$
S(g_{ij},g_{kl})
=\frac{f_{kl}f_{ji}}{f_{jl}[f_{ij},f_{kl}]}y_ig_{jl} + y_j(
\frac{f_{lj}f_{kl}f_{ji}}{f_{jl}[f_{ij},f_{kl}]}y_i -
\frac{f_{lk}f_{ij}}{[f_{ij},f_{kl}]}y_k),
$$
but $\frac{f_{lj}f_{kl}f_{ji}}{f_{jl}[f_{ij},f_{kl}]}y_i -
\frac{f_{lk}f_{ij}}{[f_{ij},f_{kl}]}y_k$ is in the first syzygy of
$f_i,f_k$, so it is a multiple of $g_{ik}$, hence
$S(g_{ij},g_{kl})$ has a standard expression with zero remainder.
Similarly, if $\frac{f_{kl}f_{ji}}{[f_{ij},f_{kl}]}$ is divided by
$f_{kl}$ or $f_{ki}$ when $i>k$, $S(g_{ij},g_{kl})$ has also a
standard expression with zero remainder.
\end{pf}

Note that $\frac{f_{kl}f_{ji}}{[f_{ij},f_{kl}]}$ is divided by
$f_{jl}$ if $[f_i,f_j]\mid f_l$. Then we have the following

\begin{Corollary}
\label{Cor} If, for any $i<j<l$, $k<l$, $k\neq i$,
$[f_{ij},f_{kl}]=1$ or $[f_i,f_j]\mid f_l$, then $f_1,\ldots,f_n$
is an s-sequence.
\end{Corollary}

\begin{Example}
\label{ex} There exist monomial sequences which do not satisfy the
sufficient condition in \cite[Proposition 1.7]{HRT}, but satisfy
our new sufficient condition. Let $f_1=x_1x_2x_3$,
$f_2=x_4x_5x_6$, $f_3=x_2x_3x_7$ and $f_4=x_7x_8x_9$. Then
$[f_{13},f_{24}]=[f_{23},f_{14}]=1$, $[f_{12},f_{34}]=x_2x_3$ but
$[f_1,f_2]\mid f_4$. Similarly, if we set $f_1=x_1^2x_3^2x_4$,
$f_2=x_1x_5^3$, $f_3=x_1^2x_4^2x_5^2$ and $f_4=x_1^2x_3x_4x_5^2$,
then $[f_{ij},f_{kl}]=1$ for all $i<j$, $k<l$, $j<l$ and $k\neq i$
but $[f_{12},f_{34}]=x_4$, and $[f_1,f_2]\mid f_4$. Hence, by
\ref{Cor}, $x_1x_2x_3, x_4x_5x_6, x_2x_3x_7, x_7x_8x_9$ and
$x_1^2x_3^2x_4, x_1x_5^3, x_1^2x_4^2x_5^2, x_1^2x_3x_4x_5^2$ are
s-sequences which do not satisfy the sufficient condition in
\cite[Proposition 1.7]{HRT}. Furthermore, the sequence
$f_1=x_2x_3x_7$, $f_2=x_1x_2x_3$, $f_3=x_3x_4x_5x_6x_7$,
$f_4=x_7x_8x_9$ satisfy all the possibilities of the condition in
\ref{s-seq} since $[f_{12},f_{34}]=1$, $[f_{13},f_{24}]=x_2$,
$[f_1,f_3]\nmid f_4$ but $f_{34}[f_{13},f_{24}]\mid f_{24}f_{31}$,
and $[f_{23},f_{14}]=x_2$, $[f_2,f_3]\nmid f_4$ but
$f_{12}[f_{23},f_{14}]\mid f_{14}f_{32}$.
\end{Example}

Recall that $f_1,\ldots,f_n$ is a strong s-sequence if it is an
s-sequence and $I_2\subseteq\cdots\subseteq I_n$, where
$I_i=(f_1,\ldots,f_{i-1}):f_i$. Note that
$I_i=(f_{1i},\ldots,f_{i-1,i})$.

\begin{Proposition}
\label{Prop}
Let $f_1,\ldots,f_n$ be a monomial sequence.
\begin{itemize}
\item[(1)] For any $i<j<k$,
$[f_i,f_j]\mid f_k$ if and only if $[f_{ik},f_{jk}]=1$.
\item[(2)] If $[f_i,f_j]\mid f_k$ for all $i<j<k$,
then $f_{1i},\ldots,f_{i-1,i}$ is a regular sequence for
$i=2,\ldots,n$.
\item[(3)] If $[f_i,f_j]\mid f_k$ for all $i<j<k$, then
$I_2\subseteq\cdots\subseteq I_n$.
\end{itemize}
\end{Proposition}

\begin{pf} (1) If $[f_i,f_j]\mid f_k$, then $[f_i,f_j]\mid
[f_i,f_k]$ and $[f_i,f_j]\mid [f_j,f_k]$, but
$f_{ik}=\frac{f_i}{[f_i,f_k]}$ and $f_{jk}=\frac{f_j}{[f_j,f_k]}$,
we see that $[f_{ik},f_{jk}]=1$. Conversely, suppose that
$[f_{ik},f_{jk}]=1$. Let $[f_i,f_j]=x_{a_1}^{\alpha_1}\cdots
x_{a_r}^{\alpha_r}$. Then, for any $1\leq t\leq r$,
$x_{a_t}^{\alpha_t}\mid [f_i,f_k]$ or $x_{a_t}^{\alpha_t}\mid
[f_j,f_k]$, hence $x_{a_t}^{\alpha_t}\mid f_k$ in two cases, it
follows that $[f_i,f_j]\mid f_k$.

(2) By (1), for any $k<l<i$, $[f_{ki},f_{li}]=1$, thus
$f_{1i},\ldots,f_{i-1,i}$ is a regular sequence.

(3) For any $i<j<k$, by assumption, $[f_i,f_j]\mid[f_i,f_k]$,
hence $f_{ik}\mid f_{ij}$. Then $I_j\subseteq I_k$, the result
follows.
\end{pf}

By \ref{Cor} and \ref{Prop}(3), we have

\begin{Corollary}
If $[f_i,f_j]\mid f_l$ for all $i<j<l$, then $f_1,\ldots,f_n$ is a
strong s-sequence.
\end{Corollary}

Note that, since $f_1,\ldots,f_n$ is a proper sequence if and only
if it is a strong s-sequence (see \cite{HRT}), we can also get the
above corollary from \ref{proper} observing that we need not to
assume that $f_1,\ldots,f_n$ is minimal in the sufficient part of
\ref{proper}.

For the sequence $f_1=x_1x_2x_3$, $f_2=x_4x_5x_6$,
$f_3=x_2x_3x_7$, $f_4=x_7x_8x_9$ in \ref{ex}, since
$[f_1,f_3]\nmid f_4$, $f_1,f_2,f_3,f_4$ is an s-sequence but not
strong.

The following proposition states that a squarefree minimal
monomial sequence is an s-sequence if it is `strong'.

\begin{Proposition}
Let $f_1,\ldots,f_n$ be a squarefree minimal monomial sequence. If
$I_2\subseteq\cdots\subseteq I_n$, then $[f_i,f_j]\mid f_k$ for
all $i<j<k$.
\end{Proposition}

\begin{pf}
By induction on $n$. When $n=3$, by $I_2=(f_{12})\subseteq
I_3=(f_{13},f_{23})$, we have $f_{12}\in (f_{13},f_{23})$. Then
$f_{13}\mid f_{12}$ or $f_{23}\mid f_{12}$. If $f_{13}\mid f_{12}$
then $[f_1,f_2]\mid [f_1,f_3]$, hence $[f_1,f_2]\mid f_3$, thus it
is enough to show that $f_{23}\nmid f_{12}$. Assume that
$f_{23}\mid f_{12}$. Then $f_2\mid f_{12}[f_2,f_3]$, but
$[f_2,f_{12}]=1$ as $f_1,f_2$ are squarefree, we have $f_2\mid
[f_2,f_3]$, hence $f_2\mid f_3$, a contradiction.

Now assume that $n>3$ and $[f_i,f_j]\mid f_k$ for all $i<j<k<n$.
Let us show that $[f_i,f_j]\mid f_n$ for any $i<j<n$. Note that we
may assume that $j=n-1$, since if $j<n-1$ then by $[f_i,f_j]\mid
f_{n-1}$ we obtain that $[f_i,f_j]\mid [f_i,f_{n-1}]$. By
assumption
$$
I_{n-1}=(f_{1,n-1},\ldots,f_{n-2,n-1})\subseteq
I_n=(f_{1n},\ldots,f_{n-1,n}),
$$
we have $f_{i,n-1}\in (f_{1n},\ldots,f_{n-1,n})$, hence
$f_{i,n-1}$ is divided by some $f_{jn}$, $1\leq j\leq n-1$. If
$f_{in}\mid f_{i,n-1}$, then $[f_i,f_{n-1}]\mid [f_i,f_n]$, hence
$[f_i,f_{n-1}]\mid f_n$. Thus it is enough to show that
$f_{jn}\nmid f_{i,n-1}$ for any $j\neq i, j\leq n-1$. Suppose that
$f_{jn}\mid f_{i,n-1}$ then $f_j[f_i,f_{n-1}]\mid f_i[f_j,f_n]$.
As $j\leq n-1$ and $i<n-1$, by induction hypothesis we have
$[f_i,f_j]\mid f_{n-1}$. Then $[f_i,f_j]\mid [f_i,f_{n-1}]$, thus
$f_j[f_i,f_j]\mid f_i[f_j,f_n]$, i.e., $f_j\mid f_{ij}[f_j,f_n]$,
but $[f_j,f_{ij}]=1$ as $f_i,f_j$ are squarefree, hence $f_j\mid
[f_j,f_n]$, so, $f_j\mid f_n$, a contradiction.
\end{pf}

\begin{Proposition} The condition in \ref{s-seq} is necessary
if $n=4$ and the sequence is squarefree.
\end{Proposition}

\begin{pf}
Suppose that $f_1,f_2,f_3,f_4$ is a squarefree s-sequence. Let us
show that $[f_{ij},f_{kl}]=1$ or
$\frac{f_{kl}f_{ji}}{[f_{ij},f_{kl}]}$ is divided by $f_{jl}$ or
$f_{ki}$ in case $i>k$ by $S(g_{ij},g_{kl})$ has a standard
expression with zero remainder.

First note that $[f_{ij},f_{jl}]=1$ as $f_1,f_2,f_3,f_4$ are
squarefree. Thus it is enough to consider $(i,j,k,l)$ is one of
$(1,2,3,4)$, $(1,3,2,4)$ and $(2,3,1,4)$. We need to show that
$[f_{12},f_{34}]=1$ or $f_{24}[f_{12},f_{34}]\mid f_{21}f_{34}$ in
the first case,  $[f_{13},f_{24}]=1$ or $f_{34}[f_{13},f_{24}]\mid
f_{31}f_{24}$ in the second case and $[f_{14},f_{23}]=1$ or
$f_{34}[f_{14},f_{23}]\mid f_{32}f_{14}$ or
$f_{12}[f_{14},f_{23}]\mid f_{32}f_{14}$ in the third case. We
will give the proof of the first case and omit the similar proofs
of the second and third cases. The following arguments do not
depend on the assumption that the sequence is squarefree.

Note that, to get a standard expression of $S(g_{ij},g_{kl})$ is
equivalent to find some $g_{st}$ whose initial term divides the
initial term of $S(g_{ij},g_{kl})$ and substitute a proper
multiple of $g_{st}$ such that the remaindered  polynomial has
smaller initial term, then do the same for this polynomial and so
on.

Suppose that $S(g_{12},g_{34})$ has a standard expression with
zero remainder. Then the initial term of
$$
S(g_{12},g_{34})=\frac{f_{21}f_{34}}{[f_{12},f_{34}]}y_1y_4-
\frac{f_{12}f_{43}}{[f_{12},f_{34}]}y_2y_3
$$
is divided by the initial term of $g_{34}$ or $g_{24}$ or
$g_{14}$, i.e., $\frac{f_{21}f_{34}}{[f_{12},f_{34}]}$ is divided
by $f_{34}$ or $f_{24}$ or $f_{14}$. In the first case, we have
$[f_{12},f_{34}]=1$. Now assume that the third case, then
$$
S(g_{12},g_{34})=\frac{f_{21}f_{34}f_{41}}{f_{14}[f_{12},f_{34}]}y_1^2-
\frac{f_{12}f_{43}}{[f_{12},f_{34}]}y_2y_3+\frac{f_{21}f_{34}}{f_{14}[f_{12},f_{34}]}y_1g_{14}.
$$
Thus $\frac{f_{12}f_{43}}{[f_{12},f_{34}]}y_2y_3$ should be
divided by the initial term of $g_{12}$ or $g_{13}$ or $g_{23}$,
i.e., $\frac{f_{12}f_{43}}{[f_{12},f_{34}]}$ is divided by
$f_{12}$ or $f_{13}$ or $f_{23}$. In the first case, we have
$[f_{12},f_{34}]=1$. In the second case, from the first equality
of $S(g_{12},g_{34})$, we have
$$
S(g_{12},g_{34})=y_1(\frac{f_{21}f_{34}}{[f_{12},f_{34}]}y_4-
\frac{f_{12}f_{43}f_{31}}{f_{13}[f_{12},f_{34}]}y_2)-
\frac{f_{12}f_{43}}{f_{13}[f_{12},f_{34}]}y_2g_{13}.
$$
As $\frac{f_{21}f_{34}}{[f_{12},f_{34}]}y_4-
\frac{f_{12}f_{43}f_{31}}{f_{13}[f_{12},f_{34}]}y_2$ is in the
first syzygy of $f_4,f_2$, we have that
$\frac{f_{21}f_{34}}{[f_{12},f_{34}]}$ is divided by $f_{24}$.
Finally, assume that $\frac{f_{12}f_{43}}{[f_{12},f_{34}]}$ is
divided by $f_{23}$. Then
$$
S(g_{12},g_{34})=\frac{f_{21}f_{34}f_{41}}{f_{14}[f_{12},f_{34}]}y_1^2
-\frac{f_{12}f_{43}f_{32}}{f_{23}[f_{12},f_{34}]}y_2^2-
\frac{f_{12}f_{43}}{f_{23}[f_{12},f_{34}]}y_2g_{23}+\frac{f_{21}f_{34}}{f_{14}[f_{12},f_{34}]}y_1g_{14}.
$$
Thus, we must have $f_{12}\mid
\frac{f_{12}f_{43}f_{32}}{f_{23}[f_{12},f_{34}]}$, hence
$\frac{f_{12}f_{43}f_{32}f_{21}}{f_{12}f_{23}[f_{12},f_{34}]}\in
R$, but
$\frac{f_{12}f_{43}f_{32}f_{21}}{f_{12}f_{23}[f_{12},f_{34}]}=
\frac{f_{12}f_{43}f_{31}}{f_{13}[f_{12},f_{34}]}$. Then
$f_{13}\mid \frac{f_{12}f_{43}}{[f_{12},f_{34}]}\cdot f_{31}$, but
$[f_{13},f_{31}]=1$, we have $f_{13}\mid
\frac{f_{12}f_{43}}{[f_{12},f_{34}]}$, reduce to the second case.
Hence, in any case, we have $[f_{12},f_{34}]=1$ or
$f_{24}[f_{12},f_{34}]\mid f_{21}f_{34}$.
\end{pf}

Finally, we consider the symmetric algebra of the ideals generated
by strong s-sequences.

\begin{Theorem}
Let $f_1,\ldots,f_n$ be a monomial sequence and
$I=(f_1,\ldots,f_n)$. Suppose that $[f_i,f_j]\mid f_k$ for all
$i<j<k$. Then
\begin{itemize}
\item[(1)] $\mbox{Sym}(I)$ is Cohen-Macaulay of dimension
$m+1$;
\item[(2)]
$e(\mbox{Sym}(I))=\sum_{i=2}^n\prod_{j=1}^{i-1}\mbox{deg}(f_{ji})$;
\item[(3)] When $f_1,\ldots,f_n$ have the same degree,
$$
\mbox{reg}(\mbox{Sym}(I))\leq\max_{2\leq i\leq
n}\{\sum_{j=1}^{i-1}\mbox{deg}(f_{ji})-(i-2)\}.
$$
\end{itemize}
\end{Theorem}

\begin{pf}
By assumption, $f_1,\ldots,f_n$ is a strong s-sequence. Then from
\ref{Prop}(2) we see that
$\mbox{dim}(R/I_i)=\mbox{depth}(R/I_i)=\mbox{dim}(R)-i+1$,
$i=2,\ldots,n$. Hence, by \cite[Propositions 2.4 and 2.6]{HRT}, we
have
\begin{eqnarray*}
\mbox{dim}(\mbox{Sym}(I))&=&\max_{2\leq i\leq n}\{\mbox{dim}(R/I_i)+i\}=\mbox{dim(R)}+1=m+1\\
\mbox{depth}(\mbox{Sym}(I))&\geq&\min_{2\leq i\leq n}\{\mbox{depth}(R/I_i)+i\}=\mbox{dim(R)}+1=m+1\\
e(\mbox{Sym}(I))&=&\sum_{2\leq i\leq n}e(R/I_i),
\end{eqnarray*}
and, when $f_1,\ldots,f_n$ have the same degree,
$$
\mbox{reg}(\mbox{Sym}(I))\leq \max_{2\leq i\leq
n}\{\mbox{reg}(I_i)\}.
$$
Then $\mbox{Sym}(I)$ is Cohen-Macaulay of dimension $m+1$, note
that from \cite[Corollary 5.10]{HSV} and the remark in front of it
we can also conclude  that
$\mbox{depth}(\mbox{Sym}(I))\geq\mbox{depth}(R)+1$. Since $I_i$ is
generated by the regular sequence $f_{1i},\ldots,f_{i-1,i}$, it is
well-known that
\begin{eqnarray*}
e(R/I_i)&=&\prod_{j=1}^{i-1}\mbox{deg}(f_{ji}),\\
\mbox{reg}(I_i)&\leq&\sum_{j=1}^{i-1}\mbox{deg}(f_{ji})-(i-2),
\end{eqnarray*}
hence
\begin{eqnarray*}
e(\mbox{Sym}(I))&=&\sum_{i=2}^n\prod_{j=1}^{i-1}\mbox{deg}(f_{ji}),\\
\mbox{reg}(\mbox{Sym}(I))&\leq&\max_{2\leq i\leq
n}\{\sum_{j=1}^{i-1}\mbox{deg}(f_{ji})-(i-2)\}.
\end{eqnarray*}
\end{pf}

\medskip

{\bf Acknowledgements}. This paper was carried out during the
author's visit to the Mathematics Department of Kansas University.
The author is grateful to Professors Craig Huneke and Daniel Katz
for helpful conversations and to the department for its
hospitality. The author is also grateful to the referee for
his/her good comments.

\medskip
\medskip
\medskip
\medskip
\medskip
\medskip
\noindent
Department of Mathematics\\
Suzhou University\\
Suzhou 215006\\
P.\ R.\ China\\
E-mail: zmtang@suda.edu.cn


\medskip
\medskip
\begin{thebibliography}{99}

\bibitem{E} D. Eisenbud,Commutative Algebra with a view towards
algebraic geometry, Graduate Texts Math. 150, Springer, 1995.

\bibitem{HRT} J. Herzog, G. Restuccia, Z. Tang, s-Sequences and symmetric
algebras, Manuscripta math. 104(2001) 479-501.

\bibitem{HSV} J. Herzog, A. Simis, W. Vasconcelos,
Approximation complexes and blowing-up rings, J. of Alg. 74(2)
(1982) 466-493.

\bibitem{Hu1} C. Huneke, The theory of $d$-sequences and
powers of ideals, Adv. in Math. 46(1982) 249-279.

\bibitem{Hu2} C. Huneke, Symbolic powers of prime ideals and special graded algebras, Comm. Alg. 9(1981) 339-366.


\bibitem{K} M.K\"uhl, On the symmetric algebra of an
ideal, Manuscripta math. 37(1982) 49-60.


\end{thebibliography}
\end{document}